\documentclass{article}

\usepackage[margin=1in]{geometry}
\usepackage{hyperref}

\usepackage[ruled,vlined,linesnumbered]{algorithm2e}
\usepackage{amsmath}
\usepackage{bm}
\usepackage{enumitem}
\usepackage{mathtools}
\usepackage{pgfplotstable}
\usepackage{siunitx}
\usepackage{subcaption}
\usepackage{tikz}
\usepackage{tikzscale}
\usepackage{tikz-network}
\usepackage{csvsimple}

\usetikzlibrary{external}
\usepgfplotslibrary{colorbrewer}
\pgfplotsset{cycle list/Set3-12,compat=newest}
\usetikzlibrary{pgfplots.colorbrewer}
\makeatletter

\begin{document}

\title{On the Emerging Potential of Quantum Annealing Hardware for Combinatorial Optimization}
\author{
    Byron Tasseff\\
    Los Alamos National Laboratory, Los Alamos, NM 87545 \\
    \and 
    Tameem Albash \\
    University of New Mexico, Albuquerque, NM 87131 \\
    \and
    Zachary Morrell, Marc Vuffray, Andrey Y. Lokhov, Sidhant Misra, Carleton Coffrin\thanks{corresponding author, cjc@lanl.gov} \\
    Los Alamos National Laboratory, Los Alamos, NM 87545
}
\date{}

\maketitle

\begin{abstract}
Over the past decade, the usefulness of quantum annealing hardware for combinatorial optimization has been the subject of much debate.
Thus far, experimental benchmarking studies have indicated that quantum annealing hardware does not provide an irrefutable performance gain over state-of-the-art optimization methods.
However, as this hardware continues to evolve, each new iteration brings improved performance and warrants further benchmarking.
To that end, this work conducts an optimization performance assessment of D-Wave Systems' most recent {\em Advantage Performance Update} computer, which can natively solve sparse unconstrained quadratic optimization problems with over $5{,}000$ binary decision variables and $40{,}000$ quadratic terms.
We demonstrate that classes of contrived problems exist where this quantum annealer can provide run time benefits over a collection of established classical solution methods that represent the current state-of-the-art for benchmarking quantum annealing hardware.
Although this work \emph{does not} present strong evidence of an irrefutable performance benefit for this emerging optimization technology, it does exhibit encouraging progress, signaling the potential impacts on practical optimization tasks in the future.
\end{abstract}


\section{Introduction}

In the 1990s, an optimization algorithm called quantum annealing (QA) was proposed with the aim of providing a fast heuristic for solving combinatorial optimization problems \cite{PhysRevE.58.5355,Ray1989,Fin1994,Far2001,San2002}.
At a high level, QA is an analog quantum algorithm that leverages the non-classical properties of quantum systems and continuous time evolution to minimize a discrete function.
\emph{Annealing} is the process that steers the dynamics of the quantum system into an a priori unknown minimizing variable assignment of that function.
Under suitable conditions, theoretical results have shown that QA can arrive at a global optimum of the desired function \cite{Bor1928,Kat1950,Jan2007}.
These results have motivated the study of using this algorithm for combinatorial optimization over the past thirty years.

Due to the computational difficulty of simulating quantum systems \cite{Fey1982}, the study of QA remained a theoretical pursuit until 2011, when D-Wave Systems produced a quantum hardware implementation of the QA algorithm \cite{Ber2010,Joh2010,Har2010,Joh2011}. 
This represented the first time that QA could be studied on optimization problems with more than a few dozen decision variables and spurred significant interest in developing a better understanding of the QA computing model \cite{Job_2018,Hau2020,Cro2021}.

The release of D-Wave Systems' QA hardware platform also generated expectations that this new technology would quickly outperform state-of-the-art classical methods for solving well-suited combinatorial optimization problems \cite{Far2001,Far2002,San2002}.
The initial interest from the operations research community was significant. 
However, through careful comparison with both complete search solvers \cite{McGeoch:2013:EEA:2482767.2482797,puget-ibm,1306.1202} and specialized heuristics \cite{1409.3934,Boixo2014,HFS_impl_2017,Man2016,Man2018,doi:10.1126/science.1252319,PhysRevA.92.042325,PhysRevX.8.031016}, it was determined that the available QA hardware was a far cry from state-of-the-art optimization methods.
These results tempered the excitement around the QA computing model and reduced interest from the operations research community.
Since the waning of the initial excitement around QA, QA hardware has steadily improved and now features better noise characteristics \cite{PRXQuantum.3.020317,9319535,2202.05847} and quantum computers that can solve optimization problems more than fifty times larger than what was possible in 2013 \cite{dw_advantage}.

Since 2017, we have been using the benchmarking practices of operations research to track the performance of QA hardware platforms and compare the results with established optimization algorithms \cite{coffrin2019evaluating,pang2021potential}.
In previous studies of this type, this benchmarking approach ruled out any potential performance benefit for using available QA hardware platforms in hybrid optimization algorithms and practical applications, as established algorithms outperformed or were competitive with the QA hardware in both solution quality and computation time.
However, in this work, we report that with the release of D-Wave Systems' most recent {\em Advantage Performance Update} computer in 2021, our benchmarking approach can no longer rule out a potential run time performance benefit for this hardware.
In particular, we show that there exist classes of combinatorial optimization problems where this QA hardware finds high-quality solutions around fifty times faster than a wide range of heuristic algorithms under best-case QA communication overheads and around fifteen times faster under real-world QA communication overheads.
This work thus provides compelling evidence that quantum computing technology has the potential for accelerating {\em certain} combinatorial optimization tasks. 
This represents an important and necessary first condition for demonstrating that QA hardware can have an impact on solving practical optimization problems.

Although this work demonstrates encouraging results for the QA computing model, we also emphasize that it {\em does not} provide evidence of a fundamental or irrefutable performance benefit for this technology.
Indeed, it is quite possible that dramatically different heuristic algorithms \cite{DunningEtAl2018,2021arXiv211113628M} or alternative hardware technologies \cite{McM2016,Hay2019,Mat2020,2103.08464} can reduce the run time performance benefits observed in this work.
We look forward to and encourage ongoing research into benchmarking the QA computing model, as closing the performance gap presented in this work would provide significant algorithmic insights into heuristic optimization methods, benefiting a variety of practical optimization tasks.

This work begins with a brief introduction to the types of combinatorial optimization problems that can be solved with QA hardware and the established benchmarking methodology in Section \ref{section:optimization-intro}.
It then presents a summary of the key outcomes from a large-scale benchmarking study in Section \ref{section:optimization-performance-analysis}, which required hundreds of hours of compute time.
In Section \ref{section:limits}, the paper concludes with some discussion of the limitations of our results and future opportunities for QA hardware in combinatoral optimization.
Additional details regarding the experimental design, as well as further analyses of computational results, are provided in the appendices.

\section{Quantum Annealing for Combinatorial Optimization}
\label{section:optimization-intro}

Available QA hardware is designed to perform optimization of a class of problems known as \emph{Ising models}, which have historically been used as fundamental modeling tools in statistical mechanics \cite{gallavotti2013statistical}.
Ising models are characterized by the following quadratic \emph{energy}  (or objective) function of $\mathcal{N} = \{1, 2, \dots, n\}$ discrete \emph{spin} variables, $\sigma_{i} \in \{-1, 1\}, \; \forall i \in \mathcal{N}$:
\begin{equation}
    E(\sigma) = \sum_{\mathclap{(i, j) \in \mathcal{E}}} {J}_{ij} \sigma_{i} \sigma_{j} + \sum_{\mathclap{i \in \mathcal{N}}} {h}_{i} \sigma_{i} \label{equation:ising-energy},
\end{equation}
where the parameters, ${J}_{ij}$ and ${h}_{i}$, define the quadratic and linear coefficients of this function, respectively.
The edge set, $\mathcal{E} \subseteq {\mathcal N} \times {\mathcal N}$, is used to encode a specific sparsity pattern in the Ising model, which is determined by the physical system being considered.
The optimization task of interest is to find the lowest energy configuration(s) of the Ising model, i.e.,
\begin{equation}
\begin{aligned}
    & \underset{\sigma}{\text{minimize}}
    & &  E(\sigma) \\
    & \text{subject to}
    & & \sigma_{i} \in \{-1, 1\}, \, \forall i \in \mathcal{N}.
\end{aligned}
\label{equation:ising-minimize}%
\end{equation}
At first glance, the lack of constraints and limited types of variables make this optimization task appear distant from real-world applications.
However, the optimization literature on {\em quadratic unconstrained binary optimization} (QUBO), which is equivalent to minimization of an Ising model's energy function, indicates how this model can encode a wide range of practical optimization problems \cite{Kochenberger2014,10.3389/fphy.2014.00005}.

\subsection{Foundations of Quantum Annealing}
\label{sec:qa_foundation}

The central idea of QA is to leverage the properties of quantum systems to minimize discrete-valued functions, e.g., finding optimal solutions to Problem \eqref{equation:ising-minimize}.
The mathematics of QA is comprised of two key elements: (i) leveraging quantum states to lift the minimization problem into an exponentially larger space and (ii) slowly interpolating (i.e., annealing) between an initial easy problem and the target problem to find high-quality solutions to the target problem.
The quantum lifting begins by introducing, for each spin, $\sigma_i \in\{-1,1\}$, a $2^{\lvert \cal N \rvert} \times 2^{\lvert \cal N \rvert}$ dimensional matrix, $\widehat{\sigma}_i$, expressible as a Kronecker product of ${\lvert \cal N \rvert}$ $2 \times 2$ matrices,
\begin{equation}
    \widehat{\sigma}_i = 
    \underbrace{\begin{pmatrix}
    1 & 0 \\
    0 & 1
    \end{pmatrix}
        \mathop{\otimes}
    \cdots
        \mathop{\otimes}
    \begin{pmatrix}
    1 & 0 \\
    0 & 1
    \end{pmatrix}}_\text{$1$ to $i-1$}
        \mathop{\otimes}
    \underbrace{\begin{pmatrix}
    1 & 0 \\
    0 & -1
    \end{pmatrix}}_\text{$i$} 
        \mathop{\otimes}
    \underbrace{\begin{pmatrix}
    1 & 0 \\
    0 & 1
    \end{pmatrix} 
        \mathop{\otimes}
    \cdots
        \mathop{\otimes}
    \begin{pmatrix}
    1 & 0 \\
    0 & 1
    \end{pmatrix}}_\text{$i+1$ to ${\lvert \cal N \rvert}$} \label{eq:sigma_hat}.
\end{equation}
In this lifted representation, the value of a spin, $\sigma_i$, is identified with the two possible eigenvalues, $1$ and $-1$, of the matrix $\widehat{\sigma}_i$.
The quantum counterpart of the energy function defined in Equation \eqref{equation:ising-energy} is the $2^{\lvert \cal N \rvert} \times 2^{\lvert \cal N \rvert}$ matrix obtained by substituting spins, $\sigma_{i}$, with the $\widehat{\sigma}_{i}$ matrices, defined in Equation \eqref{eq:sigma_hat}, within the algebraic expression for the energy.
That is,
\begin{equation}
    \widehat{E} = \sum_{\mathclap{(i,j) \in {\cal E}}} J_{ij} \widehat{\sigma}_i \widehat{\sigma}_j + \sum_{i \in {\cal N}} h_i \widehat{\sigma}_i  \label{eq:quantum_ising}.
\end{equation}
Notice that the eigenvalues of the matrix $\widehat{E}$ are the $2^{\lvert \cal  N \rvert}$ possible energies obtained by evaluating $E(\sigma)$ from Equation \eqref{equation:ising-energy} for all possible configurations of spins.
This implies that finding the minimum eigenvalue of $\widehat{E}$ is equivalent to solving Problem \eqref{equation:ising-minimize}.
This lifting is clearly impractical in the classical computing context, as it transforms a minimization problem over $2^{\lvert \cal N \rvert}$ configurations into computing the minimum eigenvalue of a $2^{\lvert \cal N \rvert} \times 2^{\lvert \cal N \rvert}$ matrix.
The key motivation for the QA computational approach is that it is possible to model $\widehat{E}$ with only $\lvert {\cal N} \rvert$ quantum bits (qubits), so it is feasible to compute over this exponentially large matrix. 

The annealing process in QA provides a method for steering a quantum system into the a priori unknown eigenvector that minimizes Equation \eqref{eq:quantum_ising} \cite{PhysRevE.58.5355,quant-ph-0001106}.
First, the system is initialized at an a priori \emph{known} minimizing eigenvector of a simple (``easy'') energy matrix, $\widehat{E}_0$. 
After the system has been initialized, the energy matrix is interpolated from the easy problem to the target problem slowly over time.
Specifically, the energy matrix at a point during the anneal is $\widehat{E}_a(\Gamma) = (1-\Gamma)\widehat{E}_0 + \Gamma \widehat{E}$, with $\Gamma$ varying from zero to one.
The annealing time is the physical time taken by the system to evolve from $\Gamma=0$ to $\Gamma=1$.
When the anneal is complete ($\Gamma=1$), the interactions in the quantum system are described by the target energy matrix.
For suitable starting energy matrices, $\widehat{E}_0$, and a sufficiently slow annealing time, the adiabatic theorem demonstrates that a quantum system remains at the minimal eigenvector of the interpolating matrix, $\widehat{E}_a(\Gamma)$ \cite{Bor1928,Kat1950,Jan2007}, and therefore achieves the minimum energy of the target problem.

\subsection{Quantum Annealing Hardware}

The computers developed by D-Wave Systems realize the QA computational model in hardware with more than $5{,}000$ qubits.
However, the engineering challenges of building real-world quantum computers are significant and have an impact on the previously discussed theoretical model of QA.
In particular, QA hardware is an open quantum system, meaning that it is affected by environmental noise and decoherence.
The coefficients in Equation \eqref{equation:ising-energy} are constrained to the ranges, $-4 \leq {h_{i}} \leq 4$, $-1 \leq {J_{ij}} \leq 1$, and nonzero ${J_{ij}}$ values are restricted to a specific sparse lattice structure (i.e., $\mathcal{E}^H \subseteq \mathcal{E}$), which is determined by the hardware's implementation.
(See Appendix \ref{section:dw-advantage} for details.)
The D-Wave hardware documentation also highlights five sources of deviation from ideal system operations called {\em integrated control errors}, which include background susceptibility, flux noise, digital-to-analog conversion quantization, input/output system effects, and variable scale across qubits \cite{dwave_docs}.
These implementation details impact the performance of QA hardware \cite{9465651}.
Consequently, QA hardware often does not find globally optimal solutions but instead finds near-optimal solutions, e.g., within 1\% of global optimality \cite{coffrin2019evaluating}.
All of these deviations from the ideal QA setting present notable challenges for encoding and benchmarking combinatorial optimization problems with available QA hardware platforms.

\subsection{Benchmarking Quantum Annealing Hardware}

Due to the challenges associated with mapping established optimization test cases to specific QA hardware \cite{coffrin2019evaluating}, the QA benchmarking community has adopted the practice of building instance generation algorithms that are tailored to specific quantum processing units (QPUs) \cite{king2015benchmarking,PhysRevA.92.042325,1701.04579,PhysRevX.6.031015,PhysRevX.8.031016,pang2021potential}.
The majority of the proposed problem generation algorithms build Ising model instances that are defined over a specific QPU's hardware graph, i.e, $(\mathcal{N}, \mathcal{E}^H$), or subsets of this graph, which are typically referred to as {\em hardware-native} problems.

In this work, we build upon an earlier class of hardware-native instances termed {\em corrupted biased ferromagnets}, or CBFMs, as proposed by \cite{pang2021potential}.
Given the QPU graph, $(\mathcal{N}, \mathcal{E}^H$), the \ref{eq:cbfm} model adopts the following distributions for hardware-native instances:
\begin{equation}
    \begin{gathered}
    P({J}_{ij} = 0) = 0, \, P({J}_{ij} = -1) = 0.625, \, P({J}_{ij} = 0.2) = 0.375, \, \forall (i, j) \in \mathcal{E}^H \\
    P({h}_i = 0) = 0.97, \, P({h}_i = -1) = 0.02, \, P({h}_i = 1) = 0.01, \, \forall i \in \mathcal{N}.
    \end{gathered}
    \label{eq:cbfm}%
    \tag{CBFM}%
\end{equation}
Benchmarking these instances on the previous generation of D-Wave's QPU architecture (i.e., the 2000Q platform using the {\sc Chimera} hardware graph) showed promising performance against state-of-the-art classical alternatives, although a clear wall-clock run time benefit was not achieved \cite{pang2021potential}.

In this work, we design a variant of the \ref{eq:cbfm} problem class called \ref{eq:cbfm-avg}, which is tailored to D-Wave's new {\em Advantage} QPU platform.
The parameter distributions are
\begin{equation}
\begin{gathered}
    P({J}_{ij} = 0) = 0.35, \, P({J}_{ij} = -1) = 0.10, \, P({J}_{ij} = 1) = 0.55, \, \forall (i, j) \in \mathcal{E}^H \\
    P({h}_i = 0) = 0.15, \, P({h}_i = -1) = 0.85, \, P({h}_i = 1) = 0, \, \forall i \in \mathcal{N}.
\end{gathered}
\label{eq:cbfm-avg}%
\tag{CBFM-P}%
\end{equation}
The \ref{eq:cbfm-avg} parameters differ from \ref{eq:cbfm}, as the new Advantage QPU architecture features a different and denser hardware graph called {\sc Pegasus}, whose topology is detailed in Appendix \ref{section:dw-advantage}.
These parameters were discovered using a metaheuristic approach that explored different combinations of the twelve parameters in this model and sought to maximize the problem's difficulty.
In each evaluation of the metaheuristic, a combination of parameters was selected, one random instance was generated following this parameterization, and a variety of classical solution methods were executed on the instance.
The instance difficulty was determined by comparing the lower and upper bounds of solutions found by these classical solution methods.
Although this approach is naive, we found that it was sufficient for the objectives of this study.
We expect that there exist classes of more challenging hardware-native instances on the {\sc Pegasus} graph, but identifying these classes is left for future work.

\section{Optimization Performance Analysis}
\label{section:optimization-performance-analysis}

In this section, we compare the performance of the D-Wave Advantage QPU and a variety of classical algorithms for optimization of \ref{eq:cbfm-avg} Ising models.
Specifically, we consider the following established classical algorithms:
\begin{itemize}
    \item A greedy algorithm based on steepest coordinate descent (SCD) \cite{pang2021potential};
    \item An integer quadratic programming (IQP) model formulation solved using the commercial mathematical programming solver \textsc{Gurobi} \cite{Billionnet2007};
    \item Simulated annealing (SA) \cite{vanLaarhoven1987,dwave-neal};
    \item A spin-vector Monte Carlo (SVMC) algorithm, which was proposed to approximate the behavior of QA \cite{shin2014quantum};
    \item Parallel tempering with iso-energetic clustering moves (PT-ICM) \cite{Zhu2015}.
\end{itemize}
SCD and IQP are general optimization approaches, intended to serve as strawman comparisons to understand solution quality, while SA, SVMC, and PT-ICM reflect high-performance classical competitors, which provide different tradeoffs in run time and solution quality.
Details of these methods and others that were considered are discussed in Appendix \ref{section:all-solution-methods}.
All of these classical optimization algorithms were executed on a system with two Intel Xeon E5-2695 v4 processors, each with 18 cores at 2.10 GHz, and 125 GB of memory.
The parameterizations used by each algorithm in this work are also detailed in Appendix \ref{section:all-solution-methods}.

For the QA hardware comparison, we use the \texttt{Advantage\_system4.1} QPU accessed through D-Wave Systems' LEAP cloud platform.
The largest system we consider features $\lvert \mathcal{N} \rvert = 5{,}387$ discrete variables and $\lvert \mathcal{E}^{H} \rvert = 25{,}324$ quadratic coefficients in the {\sc Pegasus} topology.
Solving a hardware-native optimization problem on this platform consists of (i) programming an Ising model, (ii) repeating the annealing and read-out process a number of times, and (iii) returning the highest quality solution found over all replicates.
In this analysis, we hold the annealing time constant at $62.5$ \SIUnitSymbolMicro s, which is justified in Appendix \ref{section:annealing-time-sensitivity}.
The number of anneal-read cycles are varied between $10$ and $5{,}120$ to produce different total run times.
We also leverage the {\em spin reversal transforms} feature, provided by the LEAP platform, after every $100$ anneal-read cycles to mitigate the undesirable impacts of the aforementioned integrated control errors.
For each Ising instance, this protocol typically requires less than two seconds of QPU compute time and less than $10$ seconds of total wall-clock time.

\subsection{A Characteristic Example}
\label{subsection:a-characteristic-example}
\begin{figure}[t]
    \includegraphics[width=\textwidth]{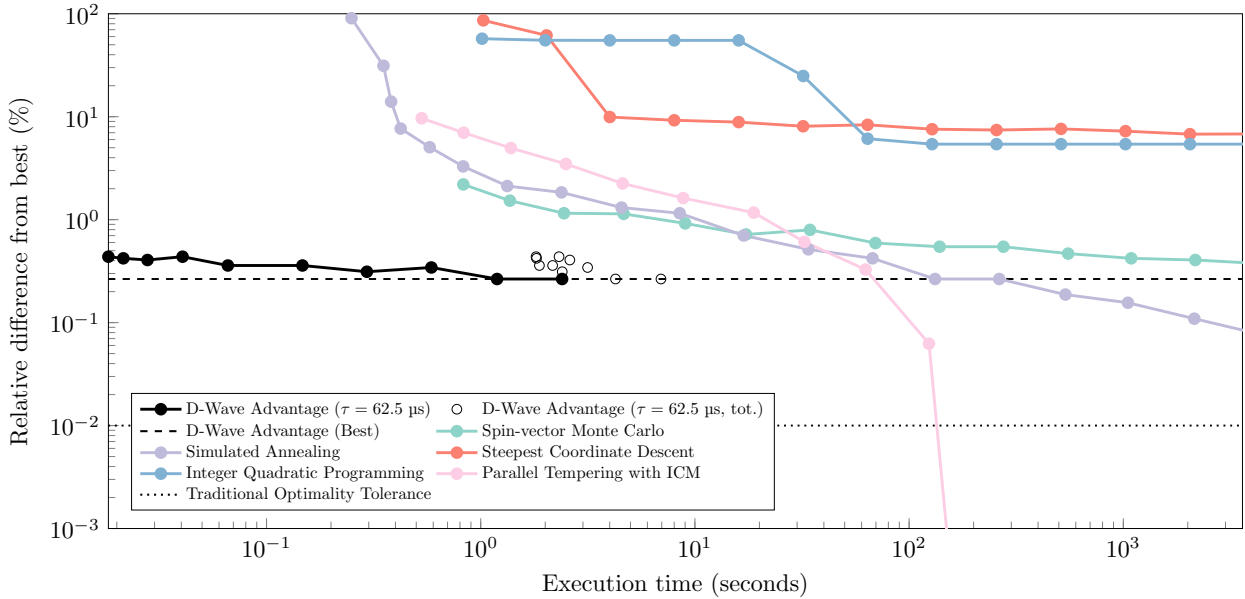}
    \caption{Evaluation of solution quality for a characteristic example of the \ref{eq:cbfm-avg} instance class with $5{,}387$ decision variables. Although the \texttt{Advantage\_system4.1} QPU does not find the best-known solution, it consistently and quickly finds solutions within $0.5\%$ of the best-known solution. Here, the dashed line corresponds to the achieved solution quality from QA when using $2{,}560$ anneal-read cycles, as used in subsequent analyses. For comparison, the dotted line corresponds to a traditional optimization tolerance of $0.01\%$, as typically used by mathematical programming solvers as a termination criterion.}
    \label{figure:characteristic-benchmark}
\end{figure}

Here, we present an evaluation of the above optimization techniques on a characteristic problem instance of the largest \ref{eq:cbfm-avg} Ising models that we considered on the \texttt{Advantage\_system4.1} QPU, with $5{,}387$ variables.\footnote{An evaluation over $50$ \ref{eq:cbfm-avg} instances in Appendix \ref{section:additional-runtime-analysis} indicates this example is not an outlier.}
For each solution technique, parameters that control the execution time of the algorithm (e.g., the number of sweeps in SA or the wall-clock time limit of the IQP method) were varied to understand their effects on solution quality.
These parameters are detailed in Appendix \ref{section:all-solution-methods}.
All other parameters remained fixed.

Benchmarking results for the \ref{eq:cbfm-avg} instance ``16'' are shown in Figure \ref{figure:characteristic-benchmark}.
Here, the horizontal axis measures the execution time of each algorithm, where each point indicates the best solution at the end of an independent algorithm execution with some set termination criterion.
The vertical axis measures the solution quality as the relative difference from the best-known solution.
Specifically, each solution's relative difference is computed as
\begin{equation}
    \textrm{\% Relative Difference} = 100\% \left(\frac{\lvert \bar{E} - E^{*} \rvert}{\lvert E^{*} \rvert}\right),
\end{equation}
where $E^{*}$ is the best-known objective value, i.e., the energy of Equation \eqref{equation:ising-energy} for the best-known solution, and $\bar{E}$ is the objective value obtained for a specific solver and execution time.

In Figure \ref{figure:characteristic-benchmark}, we first observe that the QPU (i.e., the solid black line) is shown to find high-quality, but not optimal, solutions at very fast timescales (between $0.01$ and $2.40$ seconds), with relative quality differences between $0.2\%$ and $0.5\%$ of the best-known solution.
Note that each execution time comprising the solid black line reflects a setting where the classical computer is colocated with and has exclusive access to the QPU.
In practice, QPU access is managed by D-Wave's remote cloud service, LEAP, which has overheads in both communication and job scheduling.
These solve times are reflected by the open points in Figure \ref{figure:characteristic-benchmark}, which add between one and five seconds of overhead to the total idealized solve times.
Impressively, even accounting for these significant overheads, the QPU is still able to obtain high-quality solutions well before all other classical methods that are considered.

Although the \texttt{Advantage\_system4.1} QPU is capable of quickly obtaining high-quality solutions in short amounts of time, it appears to reach a solution quality limit around $0.2\%$.
This relative difference is over an order of magnitude larger than the standard termination criterion used by mathematical programming solvers, i.e., an optimality gap of $0.01\%$ or less, delineated by the dotted line in Figure \ref{figure:characteristic-benchmark}.
To facilitate a comparison of the run time performance gained by the \texttt{Advantage\_system4.1} QPU, we thus propose a measurement that evaluates the ability of classical algorithms to \emph{match} the solution quality found by the QPU after $2{,}560$ anneal-read cycles.
The measurement we use is similar to determining the intersection with the dashed line in Figure \ref{figure:characteristic-benchmark}, albeit on a linear instead of logarithmic scale.

In this instance example, the best solution obtained by the QPU after $2{,}560$ anneal-read cycles is discovered after around $1.2$ seconds when neglecting overheads and $4.3$ seconds when including overheads.
Most solution techniques (i.e., SVMC, IQP, and SCD) do not reach this solution quality after one hour of computation.
Simulated annealing matches this quality after around $132$ seconds and, linearly interpolating between the two nearest points before and after the intersection, PT-ICM is estimated to match this quality after around $77$ seconds.
Thus, the {\em best-case} performance of the QPU in this experiment, which assumes colocation with and direct access to the QPU, provides a $64$ times improvement in run time, from $77$ seconds with PT-ICM to $1.2$ seconds.
A similar comparison using the wall-clock run time yields an improvement of around $18$ times.
That is, even when including the overhead of communicating with D-Wave's LEAP cloud service, the QPU is capable of providing a high-quality solution over an order of magnitude faster than all tested classical methods.

\subsection{Problem Scaling Run Time Trends}
\label{subsection:problem-scaling-run time-trends}
In this subsection, we investigate how the run time performance of the QPU is impacted by the size of the problem that is considered.
Unlike the previous section, here, we consider solution statistics that are \emph{aggregated} over $50$ distinct \ref{eq:cbfm-avg} instances per problem size.
Similar to the run time ratios discussed in Section \ref{subsection:a-characteristic-example}, we estimate the amount of time required for the classical algorithms to match, on average, the solution quality reported by the QPU using an annealing time of $62.5$ \SIUnitSymbolMicro s and $2{,}560$ anneal-read cycles.
This experiment is performed for {\sc Pegasus} lattice sizes ranging from two to sixteen, yielding problems with $40$ to $5{,}387$ decision variables.
For each instance, if a classical algorithm exactly matches the best solution (objective) found by the QPU, this {\em time-to-match} measurement is the earliest solve time at which that solution is obtained.
If a classical algorithm finds a solution that does not strictly match but is \emph{better} than the solution found using the QPU, the time-to-match measurement is estimated via a linear interpolation between the time at which the better solution is obtained and the time at which the worse solution, preceding it, is obtained.

\begin{figure}[t]
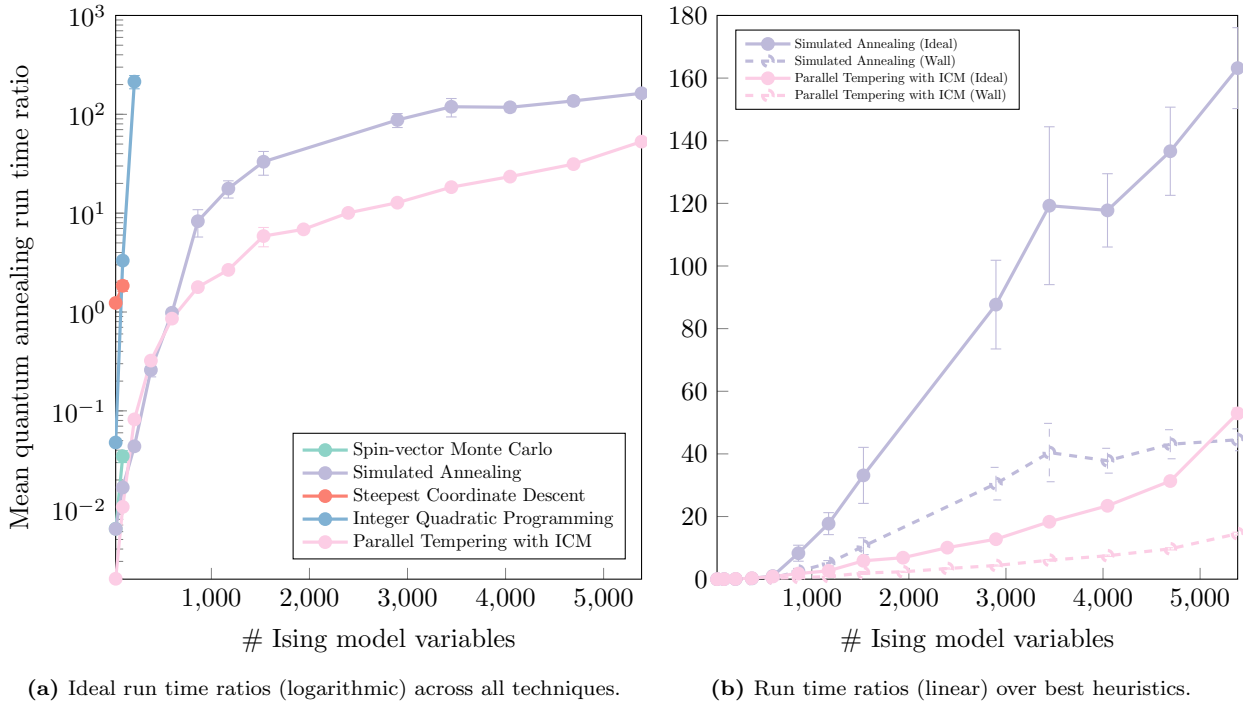

    \centering
    \begin{subfigure}[t]{0.52\linewidth}
        \centering
        \includegraphics[width=\textwidth]{figures/speedup-semilogy}
        \caption{Ideal run time ratios (logarithmic) across all techniques.}
        \label{subfigure:quantum-annealing-speedups-logarithmic}
    \end{subfigure}
    \hfill
    \begin{subfigure}[t]{0.47\linewidth}
        \centering
        \includegraphics[width=\textwidth]{figures/speedup-competitive.tikz}
        \caption{Run time ratios (linear) over best heuristics.}
        \label{figure:quantum-annealing-speedups-linear}
    \end{subfigure}
    \vspace{0.5em}
    \caption{Estimated relative computation times (``run time ratios'') required for classical algorithms to match the solution quality of the \texttt{Advantage\_system4.1} QPU as a function of problem size. Each point represents a mean computed over $50$ random \ref{eq:cbfm-avg} instances per problem size. Error bars correspond to standard errors of each mean run time ratio, and points are plotted only if the solver could match QA objectives for all $50$ instances. The computational benefits of the QPU begin to become apparent for problems with around $1{,}000$ decision variables, and these run time benefits increase steadily with problem size.}
    \label{figure:quantum-annealing-speedups}
\end{figure}

The results of conducting these scaling experiments and analyses are summarized in Figure \ref{figure:quantum-annealing-speedups}.
Figure \ref{subfigure:quantum-annealing-speedups-logarithmic} illustrates the idealized run time improvements (i.e., without including communication overheads) as a function of the number of variables in the Ising model.
It is clear that the problem scale has a large impact on the usefulness of QA hardware.
For small problem sizes (e.g., less than $500$ variables), the QPU run time is greater than most of the classical algorithms, as indicated by a run time ratio less than $10^{0}$.
For instances containing roughly $1{,}000$ variables or more, the QA hardware begins to have run time performance benefits, and only the best classical heuristics, i.e., PT-ICM and SA, are capable of matching the QA hardware's solution quality, albeit sometimes after a significant amount of time.
Note that points are excluded from Figure \ref{subfigure:quantum-annealing-speedups-logarithmic} if a solver did not match the QA solution quality for all $50$ instances.
For example, SA matched the solution quality for only $49$ of $50$ instances for \textsc{Pegasus} lattice sizes of $9$ and $10$, and hence these points are excluded from this plot.

For the two most competitive classical methods, PT-ICM and SA, Figure \ref{figure:quantum-annealing-speedups-linear} shows that the run time benefits of the QA hardware increase steadily with problem size.
This trend holds when considering both the idealized computation setting (solid lines) and the real-world setting that includes communication and scheduling overheads (dashed lines).
In particular, the estimated $15 \times$ run time ratio for the largest problem size is encouraging, as this suggests that the solutions identified by the QPU, even when accessed via a cloud computing service, can be obtained quickly enough to accelerate the performance of classical heuristic methods.

\section{Limitations and Opportunities}
\label{section:limits}
Sections \ref{subsection:a-characteristic-example} and \ref{subsection:problem-scaling-run time-trends} provide evidence that there exist classes of Ising models where available QA hardware can provide run time performance improvements over classical alternatives.
This is an encouraging result, but it is also important to recognize some limitations of this study and available QA hardware.

\paragraph{Limitations of This Study:}
The foremost limitation of this work is that it considers Ising models that are hardware-native.
Such models provide best-case scenarios for QA hardware and, thus far, have not reflected sparsity patterns of realistic combinatorial optimization tasks.
Although this work demonstrates an important \emph{necessary} condition for having a performance benefit on practical problems, it is not a \emph{sufficient} condition. Benchmarking real-world problems is required to show that these benefits can be also realized in that context.

We also note that most of the classical algorithms employed in this work did not effectively exploit parallelism, and all except SVMC and PT-ICM used their single-threaded variants.
Parallelism of classical algorithms may reduce or eliminate the performance benefits presented in this work.
Further, the benchmarks considered in this study did not evaluate other novel computing technologies or special-purpose hardware (e.g., \cite{McM2016,Hay2019,Mat2020,Hon2021,2103.08464}), which could provide improved performance on \ref{eq:cbfm-avg} instances.
Both of these avenues should be explored in future work to improve heuristic algorithms and better exploit computational resources.

Finally, we also recognize that this study does not attempt to demonstrate nor assert the much sought-after {\em scaling advantage} from quantum annealing \cite{doi:10.1126/science.1252319}, even for the contrived class of \ref{eq:cbfm-avg} instances that are considered.
This work has provided encouraging initial evidence of a class of Ising models where QA hardware can provide a practical, constant factor performance improvement over available classical algorithms.

\paragraph{Limitations of Current QA Hardware:}
The primary limitation of the QA hardware identified in this study is that it appears to approach a limit on solution quality for the largest \ref{eq:cbfm-avg} instance class, i.e., around $0.20\%$ from the best-known solution.
More evidence for this behavior is provided in Appendix \ref{section:annealing-time-sensitivity}.
As such, this work adopted a {\em time-to-match} measurement of performance, which is atypical for an optimization benchmarking study.
Additional research is required to develop extensions of the simple QA optimization protocols used in this work to understand if the hardware can achieve solutions that are within 0.01\% of global optimality, which would make this hardware's performance consistent with standard optimality tolerances used by commercial optimization tools.
QA hardware improvements to reduce noise and integrated control errors would also serve to further close this gap.

\paragraph{Future Opportunities:}
Despite the limitations of this work and current QA hardware, our results provide encouraging evidence that QA hardware is reaching a point where existing classical optimization algorithms can be practically outperformed, especially over very short timescales.
If QA hardware continues to increase in the number of qubits and hardware graph connectivity while also reducing noise properties, it is reasonable to expect that the performance gap on hardware-native problem instances will continue to increase, as suggested by the results of Section \ref{subsection:problem-scaling-run time-trends}.
Recently, D-Wave Systems announced their plans to develop an \textsc{Advantage 2} QPU, which will support over $7{,}000$ decision variables and a denser hardware graph \cite{clarity}.
If the trends observed in this work continue on this new platform, identifying even more dramatic performance gains should be possible.
Acknowledging these anticipated hardware improvements, as well as the empirical findings of this study, revisiting the topic of demonstrating a QA scaling advantage (as in \cite{PhysRevX.6.031015,PhysRevX.8.031016}) is a natural next step to establish a stronger case for a fundamental performance benefit of the QA computing model for combinatorial optimization.

\section{Conclusion}
After roughly twenty years of research and development and ten years of focused commercial development, we believe quantum annealing technology has reached a level of technical maturity and performance that warrants serious consideration by the operations research community.
This work has shown, for the first time, an order-of-magnitude run time performance benefit for quantum annealing over a wide range of classical alternatives, even when accounting for the substantial overheads involved in the practical usage of commercial quantum annealing services.
Nonetheless, significant open challenges remain in translating these performance results into benefits for practical optimization tasks.
There may be significant unrealized opportunities to hybridize this new computing technology into existing mathematical programming algorithms and impact real-world optimization challenges.
We sincerely hope that this work will inspire the operations research community to increase its consideration of the quantum annealing computing model and continue exploring how it can potentially benefit mathematical optimization algorithms and practical applications.

\section{Acknowledgment}
All work at Los Alamos National Laboratory was conducted under the auspices of the National Nuclear Security Administration of the U.S. Department of Energy under Contract No. 89233218CNA000001.
This research used resources provided by the Los Alamos National Laboratory Institutional Computing Program and was supported by the Laboratory Directed Research and Development program under the projects 20210114ER and 20210674ECR.
This material is based upon work supported by the National Science Foundation under Grant No. 2037755.

\bibliographystyle{plain}
\bibliography{bibliography}

LA-UR-22-29705

\appendix

\section{D-Wave Advantage System Details}
\label{section:dw-advantage}

\begin{figure}[t]
    \includegraphics[width=\textwidth]{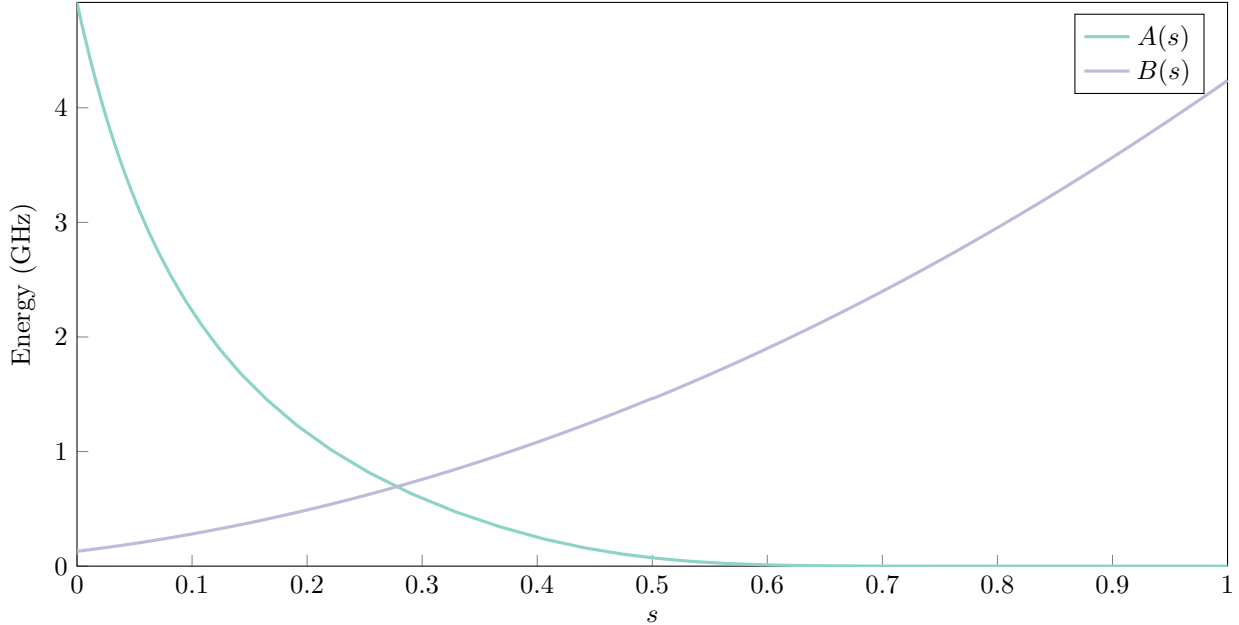}
    \caption{Annealing schedule used by D-Wave Systems' \texttt{Advantage\_system4.1} QPU used throughout this study, in units where the Planck constant set to one.}
    \label{figure:dw-annealing-schedule}
\end{figure}

\begin{figure}[t]
    \centering
    \includegraphics[width=0.9\textwidth]{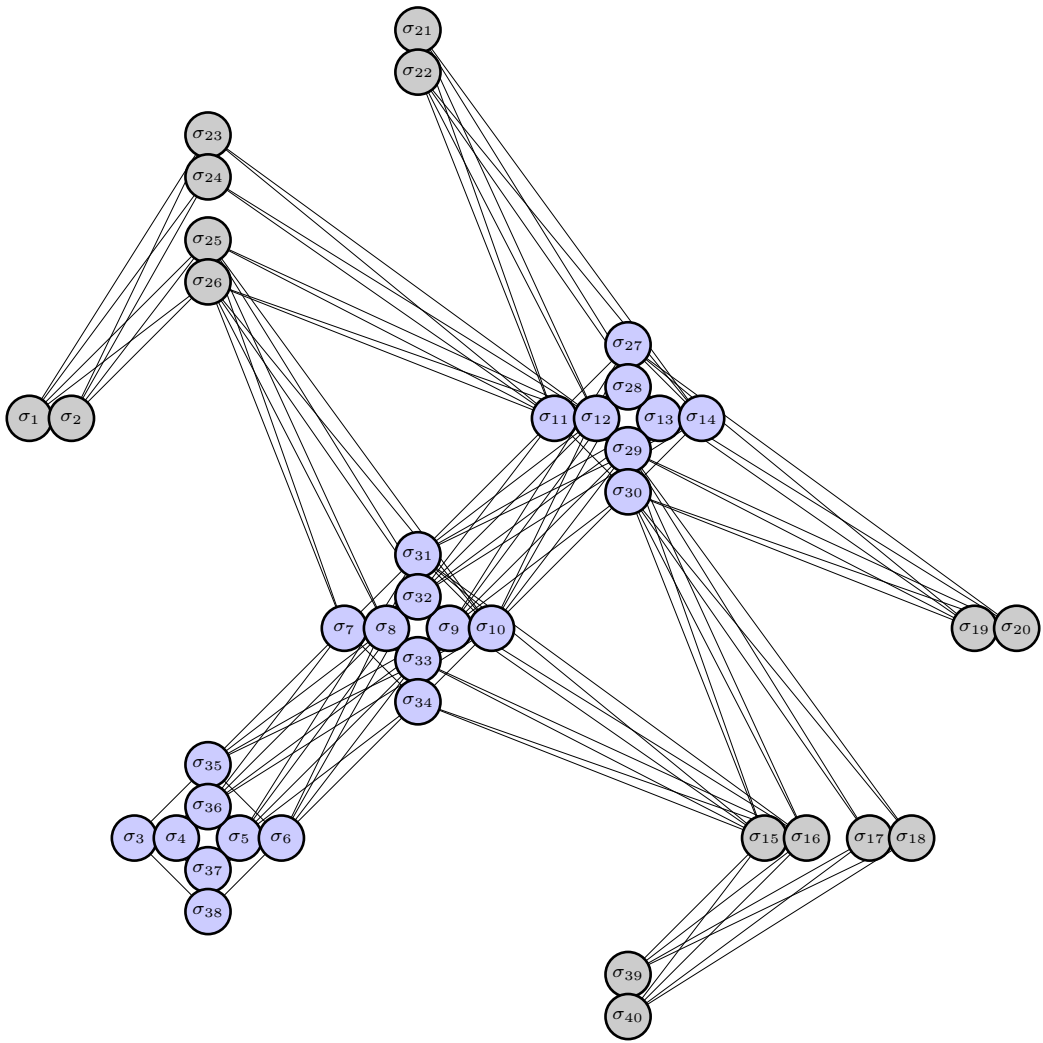}
    \caption{Topology of a \textsc{Pegasus} graph with a size parameter of two, where blue nodes comprise a unit cell.}
    \label{figure:pegasus-topology-lattice-size-2}
\end{figure}

In an ideal setting, a QA hardware device would provide the user control over the Ising models and annealing schedules that are executed.
However, available QA hardware impose restrictions on the classes of Ising models and annealing schedules that can be considered per the hardware's architecture.
In this paper, we used the \texttt{Advantage\_system4.1} QPU, which is based on the \textsc{Advantage} architecture \cite{Boo2020}.
This system implements the transverse-field Ising model,
\begin{equation}
E(s) = - A(s) \sum_{i \in \mathcal{N}} \widehat{\sigma}^{x}_i + B(s) \left( \sum_{(i, j) \in \mathcal{E}^{H}} J_{ij} \widehat{\sigma}^{z}_i \widehat{\sigma}^{z}_j - \sum_{i \in \mathcal{N}}  h_{i} \widehat{\sigma}^{z}_i \right), \label{eq:tising}%
\end{equation}
where $\widehat{\sigma}^{x}$ and $\widehat{\sigma}^{z}$ are the Pauli operators for the $x$ and $z$ bases, respectively, and the Ising model of interest is encoded in the $z$ basis. 
The hardware has a fixed global annealing schedule shown in Figure \ref{figure:dw-annealing-schedule} and a qubit connectivity graph known as {\sc Pegasus}, containing a total of 5,387 qubits.
Figure \ref{figure:pegasus-topology-lattice-size-2} illustrates a \textsc{Pegasus} graph topology with a lattice ``size parameter'' of two.
A \textsc{Pegasus} \emph{unit cell} contains twenty-four qubits (here, the three diagonal clusters of eight), with each qubit coupled to one similarly aligned qubit \emph{within} the cell and two similarly aligned qubits \emph{in adjacent cells}.
An \textsc{Advantage} QPU is a lattice of $16 \times 16$ such unit cells.
This example illustrates that the use of QA is limited to Ising models that are \emph{topologically consistent} with the QPU used and \emph{size-restricted} by its number of qubits.
In this work, we only consider Ising models that are subsets of this hardware graph, i.e., hardware-native Ising models.
This allows us to avoid effects caused by the {\em embedding} problem, which is defined as mapping a desired problem structure onto this hardware graph.

Even if the Ising model is topologically consistent, achieving a globally optimal solution to Problem \eqref{equation:ising-minimize} using QA is often difficult at large scales.
As explained in Section \ref{section:optimization-intro}, QA has been theoretically shown to minimize the energy of specific target Ising models.
Experimentally, however, this is not always the case.
One reason is the unavoidable corruption of the quantum system from its external environment.
For example, as previously stated in Section \ref{section:optimization-intro}, QPUs are known to suffer from a wide range of \emph{integrated control errors}, including susceptibility to background noise, flux noise, digital-to-analog conversion, and variable scales and coefficient biases across qubits \cite{9465651}.

It is also important to note that QPUs, such as the \textsc{Advantage}, often distinguish between the tasks of \emph{optimization} and \emph{sampling}.
The former aims at determining a solution to Problem \eqref{equation:ising-minimize}, while the latter aims at characterizing the \emph{energy landscape} of an Ising model.
In practice, since quantum annealing does not experimentally guarantee convergence to a global optimum, optimization is performed by conducting multiple \emph{anneal-read cycles} (or ``reads'') for the target Ising model.
Each read provides an assignment of the decision variables, which may or may not be the true global optimum of Problem \eqref{equation:ising-minimize} but which tends to be high-quality (e.g., within one percent of optimality).
The best encountered solution is then extracted by computing the minimum of these multiple reads.

In this paper, we consider only two key parameters for the configuration of a QPU, and we ignore any environmental factors that are beyond our control.
The first parameter is the number of anneal-read cycles performed to obtain a solution to Problem \eqref{equation:ising-minimize}.
Here, a larger number generally corresponds to a higher likelihood of sampling the true ground state, i.e., as the number of reads increases, solution quality also tends to increase.
The second parameter is the \emph{annealing time} that defines \emph{each} anneal-read cycle, which is on the order of microseconds.
Here, longer anneal times typically encourage convergence of the QPU to higher-quality solutions, as would be expected from quantum annealing theory.
Empirical evidence of this tendency is provided in Appendix \ref{section:annealing-time-sensitivity}.

\section{Optimization Methods}
\label{section:all-solution-methods}

In the interest of clarity, Section \ref{section:optimization-performance-analysis} briefly introduced five established and distinct solution approaches for solving the proposed \ref{eq:cbfm-avg} problem instances.
However, in the development of this work, we considered ten different solution approaches.
In this section, we discuss all ten methods that were considered and provide additional details regarding how each was configured.
The algorithms are divided into two broad categories: {\em complete search methods}, which, given sufficient time, provide proofs of global optimally, and {\em local search methods}, which are heuristics with no guarantees nor bounds on solution quality.
Unless noted otherwise, each of these solvers was implemented in C/C\texttt{++}.

\subsection{Complete Search Methods}
In all of the complete search methods that are considered, we leverage the bijection of Ising models and QUBO models to convert a given Ising model of interest to its equivalent QUBO form prior to its solution \cite{coffrin2019evaluating}.
Preliminary experiments indicated that the QUBO formulation is advantageous for commercial mathematical programming tools, which include specialized techniques for solving problems with discrete $\{0, 1\}$ decision variables.
In this work, we considered both \textsc{Gurobi}, version 9.1, and \textsc{CPLEX}, version 22.1, for solving these mathematical programming models.
Both solvers were parameterized to use only one thread during solution.
Solver execution time limits were varied from one to $8{,}192$ seconds in multiples of two for each instance.
All other solver parameters assumed their default values.
As discussed by \cite{puget-ibm}, we acknowledge that better performance may be attained with a larger number of threads and adjusting the solvers' other parameter settings.

Mathematical programming solution approaches were considered by some of the earliest QA benchmarking studies \cite{McGeoch:2013:EEA:2482767.2482797} and are widely recognized as an unfair point of comparison, as specialized heuristics can be much faster than these general purpose solvers \cite{puget-ibm}.
However, it was shown by \cite{coffrin2019evaluating} and \cite{pang2021potential} that these solvers provide an important point of reference for ruling out easy optimization tasks, as is similarly done in this study.

\subsubsection{Integer quadratic programming}
\label{subsubsection:integer-quadratic-programming}
The integer quadratic programming (IQP) method considered in this study consists of using black-box commercial mathematical programming software (in this case, \textsc{Gurobi} and \textsc{CPLEX}) to pose and solve the following QUBO Ising model translation:
\begin{equation}
\begin{aligned}
    & \text{minimize}
    & & \sum_{\mathclap{(i, j) \in \mathcal{E}}} {c}_{ij} x_{i} x_{j} + \sum_{\mathclap{i \in \mathcal{N}}} {c}_{i} x_{i} + {c} \\
    & & & x_{i} \in \{0, 1\}, \, \forall i \in \mathcal{N},
\end{aligned}%
\label{equation:qubo-minimize}%
\end{equation}
where $c_{ij}$, $c_{i}$, and $c$ are QUBO parameters derived from the Ising model parameters of Equation \eqref{equation:ising-energy}.

\subsubsection{Integer linear programming}
\label{subsubsection:integer-linear-programming}
The integer linear programming (ILP) model is an equivalent reformulation of the IQP model, where each variable product, $x_{i} x_{j}$, is lifted into a new variable space, $x_{ij}$.
For each new variable, constraints are used to model each conjunctive relationship, i.e., $x_{ij} = x_i \land x_j$.
Notably, this new problem contains no quadratic terms in the objective, and all constraints are affine.
The ILP reformulation of the IQP formulation discussed in Section \ref{subsubsection:integer-quadratic-programming} is
\begin{equation}
\begin{aligned}
    & \text{minimize}
    & & \sum_{\mathclap{(i, j) \in \mathcal{E}}} {c}_{ij} x_{ij} + \sum_{\mathclap{i \in \mathcal{N}}} {c}_{i} x_{i} + {c} \\
    & \text{subject to}
    & & x_{i} + x_{j} - 1 \leq x_{ij}, \, \forall (i, j) \in \mathcal{E} \\
    & & & x_{ij} \leq x_{i}, \, \forall (i, j) \in \mathcal{E} \\
    & & & x_{ij} \leq x_{j}, \, \forall (i, j) \in \mathcal{E} \\
    & & & x_{i} \in \{0, 1\}, \, \forall i \in \mathcal{N} \\
    & & & x_{ij} \in \{0, 1\}, \, \forall (i, j) \in \mathcal{E}.
\end{aligned}%
\label{equation:ilp-formulation}%
\end{equation}
This model was also evaluated in some of the earliest QA and QUBO benchmarking studies \cite{puget-ibm,1306.1202,Billionnet2007}, which suggested that it is a very effective formulation for solving Ising models with sparse graphs (e.g., D-Wave's previous {\sc Chimera} topology).
Like \cite{pang2021potential}, our results indicate that IQP solvers have improved significantly over the past decade and consistently outperform ILP methods for solving QUBO problems on large problem instances.
We remark that, in this study, only \textsc{Gurobi} (i.e., not \textsc{CPLEX}) was considered when solving ILP reformulations.
Substantial performance differences were not anticipated when solving with \textsc{CPLEX}.

\subsection{Local Search Methods}
Although global optimization methods (e.g., IQP and ILP) are useful for measuring the quality of solutions produced by QA hardware \cite{Billionnet2007,coffrin2019evaluating,Baccari_2020}, it is broadly accepted that local search algorithms are the most appropriate points of computational comparison for QA methods \cite{aaronson-insert-d-wave-post}.
Given that an enumeration of all local search methods would be an incredibly large and impractical undertaking, this work focuses on a set of the most common or readily available approaches: greedy search via steepest coordinate descent, tabu search, messaging passing, Markov chain Monte Carlo, simulated annealing, spin-vector Monte Carlo, and parallel tempering with isoenergetic cluster moves.
Next, each method is briefly described.

\subsubsection{Steepest coordinate descent}
\label{subsubsection:steepest-coordinate-descent}
The simplest heuristic algorithm considered in this work is a steepest coordinate descent, or greedy, approach that was implemented in Julia.
This algorithm assigns variable values one-by-one, always taking an assignment that most decreases the objective.
Specifically, the approach begins with unassigned values, i.e., $\sigma_{i} = 0$, for all $i \in \mathcal{N}$, then repeatedly applies the following assignment logic until all variables have been assigned a value of $-1$ or $1$:
\begin{subequations}
\begin{align}
    i, v &= \arg\min\left\{E(\sigma_1, \ldots, \sigma_{i-1}, v, \sigma_{i+1}, \ldots,\sigma_N) : i \in \mathcal{N}, \, v \in \{-1, 1\}\right\} \\
    \sigma_i &= v
\end{align}%
\end{subequations}
In each application of the above rule, ties in the $\arg\min$ are broken at random, giving rise to potentially stochastic outcomes.
Once all variables have been assigned their values via this rule, the algorithm is repeated until a run time limit has been reached, and only the best discovered solution is returned.
This approach is fast and effective on Ising models with minimal amounts of frustration.
In our work, the wall-clock time limit of the steepest coordinate descent algorithm was varied from one to $8{,}192$ seconds in multiples of two for each problem instance to understand solution improvement.

\subsubsection{Tabu search}
Tabu search is a metaheuristic commonly used for solving combinatorial optimization problems \cite{Glover1998}.
Like other local methods, tabu search considers a ``move'' operator that, for a given solution, generates a number of other possible solutions (i.e., the solution's ``neighborhood'').
The ``best'' solution from this neighborhood is then selected, and the process repeats.
To prevent cycling, a list of ``tabu moves'' is employed, which prohibits moves that would lead to revisiting a previously-encountered solution.
This list is updated as the algorithm progresses, i.e., new moves are added to the list, and moves with a large ``tabu tenure'' are discarded.
Tabu moves are, however, sometimes allowed if they lead to an improved solution \cite{beasley1998heuristic}.

In this study, we used an open-source implementation of the MST2 multistart tabu search algorithm of \cite{palubeckis_multistart_2004}, provided by \cite{dwave-tabu}, which is specialized for solving QUBOs.
For brevity, we forgo details of this algorithm.
However, in our work, for each problem instance, the number of ``reads,'' \texttt{num\_reads}, was varied between one and $2{,}048$ in multiples of two, and the total run time per ``read,'' \texttt{timeout}, was limited to a fixed value of one second.
For the largest instances studied, the most time-intensive parameterizations resulted in execution times of nearly one hour.

\subsubsection{Message passing}
The message-based min-sum (MS) algorithm is an adaptation of the belief propagation algorithm for solving minimization problems over networks \cite{768759,mezard2009information}.
A key property of the MS algorithm is its ability to identify global minima of cost functions over networks with tree-like structures, i.e., if no cycles are formed by the interactions in $\mathcal{E}$.
In the more general case, MS is a heuristic minimization method \cite{mezard2009information}.
Nonetheless, it remains a popular heuristic technique, favored in communication systems for its low computational cost and performance on tree-like networks \cite{vuffray2014cavity}.

For the model we consider, i.e., Problem \eqref{equation:ising-minimize}, the MS messages, $\epsilon_{i \rightarrow j}$,  are computed iteratively along \emph{directed} edges, $i \rightarrow j$ and $j \rightarrow i$, for each $(i, j) \in \mathcal{E}$ according to the min-sum equations
\begin{subequations}
\begin{align}
    \mathrm{SSL}(x, y) &= \min(x, y)-\min(-x, y) - x \\
    \epsilon^{t + 1}_{i \rightarrow j} &= \mathrm{SSL}\left(2 {J}_{ij}, 2 {h}_i +
    \sum_{\mathclap{k \in  \mathcal{E}(i) \setminus j}} \epsilon^{t}_{k \rightarrow j}\right)
    \label{equation:min_sum)}.
\end{align}
\end{subequations}
Here, $\mathcal{E}(i) \setminus j$ denotes the neighbors of $i$ without $j$, and $\mathrm{SSL}$ denotes the symmetric saturating linear (SSL) transfer function.
Once a fixed point of Equation \eqref{equation:min_sum)} is obtained or a prescribed termination criterion is achieved (e.g., a time limit), the MS algorithm outputs the configuration
\begin{equation}
    \bar{\sigma}_{i} = - \mathrm{sign} \left(2 {h}_i + \sum_{k \in \mathcal{E}(i)} \epsilon_{k \rightarrow j} \right) \label{equation:min_sum_assignment}.
\end{equation}
If the argument of $\mathrm{sign}$ is zero, a value of $1$ or $-1$ is assigned randomly with equal probability.

This algorithm is implemented in the Python programming language.
In our work, the wall-clock time limit of the algorithm was varied from one to $8{,}192$ seconds in multiples of two for each instance.

\subsubsection{Markov chain Monte Carlo}
Markov chain Monte Carlo (MCMC) algorithms are a class of methods to generate samples from complex probability distributions.
A natural MCMC method for solving Ising models is given by Glauber dynamics (GD), where the value of each variable is updated according to its conditional probability distribution.
Glauber dynamics is often used as a method for producing samples from Ising models at \emph{finite temperatures} \cite{doi:10.1063/1.1703954}.
This work considers the so-called \emph{zero temperature} GD algorithm, which is the optimization variant of the GD sampling method.
This method is also used in physics as a simple model for describing avalanche phenomena in magnetic materials \cite{Dhar_1997}.
From an optimization perspective, the approach is similar to the single-variable greedy local search algorithm previously described in Appendix \ref{subsubsection:steepest-coordinate-descent}.

A step $t$ of the GD algorithm comprises selecting each spin variable, for $i \in \mathcal{N}$, in a random order and comparing the objective of the current configuration, $\bar{\sigma}^{t}$, to a configuration where $\bar{\sigma}^{t}_{i}$ is flipped in sign.
If the objective value is smaller in the flipped configuration, i.e.,
$E(\bar{\sigma}^{t}) > E(\bar{\sigma}^{t}_1, \ldots, -\bar{\sigma}^{t}_i, \ldots, \bar{\sigma}^{t}_{N})$,
then the flipped configuration is selected as the new configuration, i.e., $\bar{\sigma}^{t + 1} = (\bar{\sigma}^{t}_1, \ldots, -\bar{\sigma}^{t}_{i}, \ldots, \bar{\sigma}^{t}_{N})$.
If, after visiting all of the variables, no single-variable flip can improve the current assignment, then the configuration is identified as a local minimum.
The algorithm is then restarted with a new, randomly generated initial configuration.
This process is repeated until a run time limit is reached.

This algorithm is implemented in the Python programming language.
In our work, the wall-clock time limit of the algorithm was varied from one to $8{,}192$ seconds in multiples of two for each instance.

\subsubsection{Simulated annealing}
Simulated annealing is an algorithm inspired by the process of annealing materials to produce improved structural properties, applied to the computational task of solving a combinatorial optimization problem.
In condensed matter physics, \emph{annealing} is a physical process by which the temperature of a solid immersed in a \emph{heat bath} is first increased, and particles of the solid randomly rearrange themselves into a liquid phase.
The temperature of the heat bath is then slowly decreased, and, in an ideal case, all particles rearrange themselves to a corresponding solid lattice structure \cite{vanLaarhoven1987}.
Analogously, in simulated annealing, solutions are probabilistically perturbed according to a ``temperature'' schedule until reaching some minimum.

The algorithm originated with \cite{doi:10.1063/1.1699114}, who proposed a Monte Carlo method for simulating the evolution of a solid to thermal equilibrium at a fixed temperature, $T$.
Specifically, given the current state of the system, $\bar{\sigma}$, a small perturbation is selected.
If the difference in energy, $\Delta E$, between the current and perturbed states is negative, the process is continued with the new state.
Otherwise, the \emph{probability of acceptance} of the new state is $\exp\left({-\frac{\Delta E}{k_{B} T}}\right)$, called the Metropolis criterion, where $k_{B}$ is the Boltzmann constant.
Simulated annealing can thus be thought of as a \emph{sequence of Metropolis algorithms}, where the temperature is gradually decreased until reaching some termination criterion (e.g., a predefined temperature value) \cite{vanLaarhoven1987}.

In this work, we employ the \textsc{dwave-neal} Python package, which provides an implementation of simulated annealing for general Ising model graphs.
In each execution of the algorithm, the number of anneal-read cycles was set to one hundred, and the number of simulated annealing ``sweeps,'' or random perturbations per discrete temperature value, varied between one and $262{,}144$ in multiples of two.
All other parameters assumed their default values.
For the largest Ising models considered in this study (i.e., a \textsc{Pegasus} graph with a size parameter of sixteen), the longest wall-clock execution times when employing the most time-intensive parameterization were a little over an hour.

\subsubsection{Spin-vector Monte Carlo}
Spin-vector Monte Carlo \cite{shin2014quantum} is a classical emulator of a noisy quantum annealer, where the system of qubits, i.e., $\mathcal{N}$, is replaced by a system of two-dimensional rotors.
The system's energy is given by
\begin{equation}
    E(\theta,s) = - A(s) \sum_{i \in \mathcal{N}} \sin \theta_i + B(s) \left( \sum_{i \in \mathcal{N}} {h}_i \cos \theta_i + \sum_{(i, j) \in \mathcal{E}} {J}_{ij} \cos \theta_i \cos \theta_j \right),
\end{equation}
where $s \in [0,1]$ is the dimensionless annealing parameter; the functions $A(s)$ and $B(s)$ are taken to match those of the quantum annealer, illustrated in Figure \ref{figure:dw-annealing-schedule}; and $\left\{ \theta_i \in [0, \pi) \right\}$ are the orientations of the rotors.
In principle, the rotors can point along any orientation in $[0,2 \pi)$, but with this energy function, the energy is minimized by an orientation in $[0, \pi)$, which allows us to restrict each $\theta_{i}$ to this range of angles.
Using this energy, the rotors' orientations are updated using a fixed-temperature Metropolis-Hastings criterion.
This description can be physically motivated as a description of superconducting flux qubits in the strong system-bath interaction limit \cite{Cro2016} when using the spin-coherent path integral formalism \cite{Kla1979,Alb2015}.

The SVMC algorithm proceeds as follows.
The rotors are initialized to point along $\theta_i = \pi/2$, corresponding to the uniform superposition state.
The anneal is discretized in steps of $\Delta s$ such that $s_k = k \Delta s$, and at every value of $s_k$, a certain number of sweeps of rotor updates are performed.
A sweep corresponds to attempting one Metropolis-Hasting update for each rotor.
The update goes as follows.
For a rotor with orientation $\theta_i$, we randomly choose a new orientation $\theta_i'$, and we calculate the change in energy, $\Delta E$, associated with this new orientation and accept the update according to the Metropolis-Hastings probability \cite{doi:10.1063/1.1699114,Has1970},
\begin{equation}
   p = \min \left( 1, \exp(- \beta (E(\theta',s)- E(\theta,s))) \right),
\end{equation}
where $\beta = 3.9983$ GHz$^{-1}$ is the fixed inverse temperature of the simulation, corresponding to a physical temperature of $12$ mK in units where the Planck constant set to one.
At the end of the anneal, the rotors are projected to Ising spins, and we can measure their energy according to Equation~\eqref{equation:ising-energy}.
If $\cos \theta_i > 0$, we assign a spin value of $\sigma_i = 1$, and if $\cos \theta_i < 0$, we assign a spin value of $\sigma_i = -1$. If $\cos \theta_i$ is zero, a value of $1$ or $-1$ is assigned randomly with equal probability.

For each problem instance, we executed eight independent SVMC simulations in parallel on eight cores, reporting the lowest energy found by the eight independent runs.
For each independent SVMC run in the parallel set, the number of sweeps used by the algorithm was varied from $1{,}000$ to $8{,}192{,}000$ in multiples of two.
For the largest Ising models considered in this study, the longest wall-clock execution times when employing the most time-intensive parameterization were a little over an hour.

\subsubsection{Parallel tempering with isoenergetic cluster moves}
Parallel tempering (PT) \cite{Swe1986,Gey1991,Huk1996} is a method that uses multiple Markov chain Monte Carlo simulations to improve the equilibration dynamics.
Each simulation, known as a ``replica,'' has a unique inverse temperature, $\beta$, and its spin configuration is evolved using single-spin Monte Carlo updates.
That is, a spin with value $\sigma_i$ is flipped to $\sigma_i' = -\sigma_i$, and the update is accepted according to the Metropolis-Hastings probability \cite{doi:10.1063/1.1699114,Has1970},
\begin{equation}
   p = \min \left( 1, \exp(- \beta (E(\sigma') - E(\sigma))) \right).
\end{equation}
After each replica performs a fixed number of sweeps, where a sweep corresponds to performing a single Monte Carlo update on all the spins, a parallel tempering update is performed.
In a parallel tempering update, the spin configurations (or equivalently the temperatures) of two replicas are swapped according to the Metropolis-Hastings probability (\cite{doi:10.1063/1.1699114,Has1970}),
\begin{equation}
   p = \min \left( 1, \exp( (\beta_r - \beta_{r'}) (E(\sigma_r) - E(\sigma_{r'}))) \right),
\end{equation}
where $\beta_{r}$ and $E(\sigma_{r})$ denote the inverse temperature and energy of the spin configuration of replica $r$, respectively.
This choice preserves detailed balance.

The isoenergetic cluster move (ICM) \cite{Zhu2015} is based on the cluster update by Houdayer \cite{Hou2001}.
Two independent parallel tempering simulations are performed. After a fixed number of Monte Carlo sweeps and the replica swaps are performed, we perform the ICM update.
This move is performed on replicas with a temperature less than one.
For equal temperature replicas $r$ and $r'$ from each PT simulation, we calculate the site overlap,
\begin{equation}
    q_i = \sigma_i^{(r)} \sigma_i^{(r')}, \, i \in \mathcal{N},
\end{equation}
where $\sigma_i^{(r)}$ is the spin configuration of the $i$-th site in the $r$-th replica.
We randomly pick a spin $i$ with $q_i = -1$, and starting from this spin, we build the largest connected (according to $\mathcal{E}$) cluster of spins with $q_i = -1$.
All spins in this cluster are then flipped in both replicas, $r$ and $r'$.
According to \cite{Zhu2015}, the ICM move need not be performed on all replicas, and we only implement it on replicas with an inverse temperature greater than one.

Therefore, the complete parallel tempering with isoenergetic cluster moves (PT-ICM) algorithm proceeds as follows.
All the replicas of the two independent PT simulations are initialized with random spin configurations.
Each replica is evolved independently, with two Monte Carlo sweeps.
This is followed by a parallel tempering update for all pairs of neighboring temperatures, where the order of the parallel tempering pairs is chosen randomly.
Finally, an ICM update is performed.

\begin{figure}[t]
    \includegraphics[width=\textwidth]{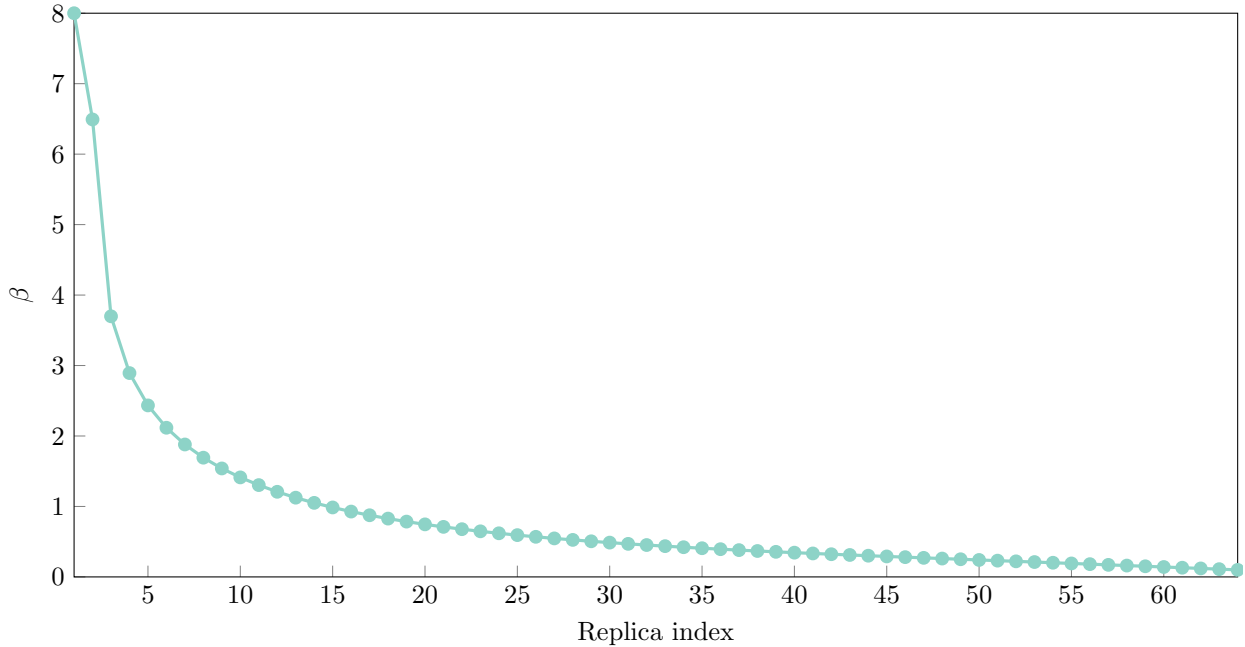}
    \caption{Inverse temperature distribution over the $64$ PT-ICM replicates considered in the ``tuned'' variant.}\label{figure:pt-beta-distribution}
\end{figure}

In our simulations, we have chosen to use 64 replicas, and we considered two inverse temperature distributions.
The first is a typical geometrical distribution, untuned with respect to the instances considered throughout this work.
The PT-ICM algorithm using this ``default'' inverse temperature distribution is represented in Figures \ref{figure:characteristic-benchmark} and \ref{figure:quantum-annealing-speedups}, and it represents an untuned variant of the algorithm, similar to the other methods that we compare.
The second ``non-default'' distribution considers values ranging between $[0.1, 8]$ to achieve a swap probability of approximately 0.23 \cite{Kon2005} between most neighboring temperature replicas at the largest system size we study, using the scheme of \cite{Roz2019}.
We illustrate this distribution in Figure \ref{figure:pt-beta-distribution}.
Further, PT-ICM results when using this distribution are postfixed with (Opt.) in subsequent figures of this appendix.
We remark, however, that this choice of temperatures is not necessarily optimal \cite{Kat2006}, but it is generally considered a reasonable choice.
In principle, at the smaller system sizes, we require fewer replicas to achieve the same swap probability, but we do not perform this optimization, here.

In our reported results, for each problem instance, we executed eight independent PT-ICM simulations in parallel on eight cores, and we reported the lowest energy found by the eight independent runs.
For each independent PT-ICM run in the parallel set, the number of parallel tempering updates used by the algorithm, prior to termination, was varied from two to $16{,}384$ in multiples of two.
For the largest Ising models considered in this study, the longest wall-clock execution times when employing the most time-intensive parameterization were typically less than thirty minutes.

\section{Annealing Time Sensitivity}
\label{section:annealing-time-sensitivity}

In the quantum annealing solution approach, the annealing time, $\tau$ (\SIUnitSymbolMicro s), of the algorithm  plays a critical role in solution quality \cite{doi:10.1126/science.1252319,PhysRevX.8.031016}, since it is generally the case that increasing the annealing time will increase quality until a globally optimal solution is achieved.
Consequently, in this work, we considered different annealing times to understand their impacts on the quality of \ref{eq:cbfm-avg} solutions.
This also allowed us to select a suitable value of $\tau$.

\begin{figure}[t]
    \includegraphics[width=\textwidth]{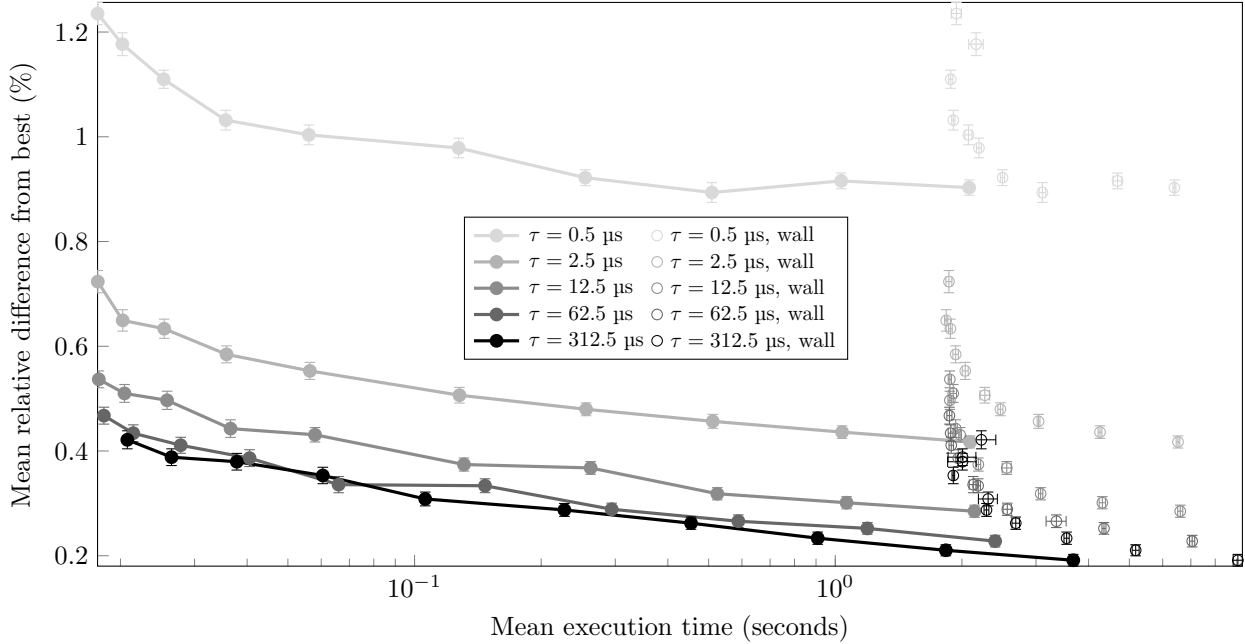}
    \caption{Effects of annealing time on solution quality and ideal/wall-clock execution time for QA hardware, averaged over $50$ \ref{eq:cbfm-avg} instances with $5{,}387$ decision variables ({\sc Pegasus} lattice size sixteen). Here, each point corresponds to a mean relative difference, computed over $50$ individual relative differences from the best solution obtained for each instance. Vertical error bars correspond to the standard error of each mean. Horizontal error bars correspond to standard errors of each mean solve time, each of which is based on a particular QPU execution time parameterization (e.g., $2{,}560$ anneal-read cycles).}
    \label{figure:annealing-time-performance}
\end{figure}

Figure \ref{figure:annealing-time-performance} shows the mean relative performance of the QA hardware over a set of $50$ instances with $5{,}387$ decision variables ({\sc Pegasus} lattice size sixteen).
Here, the number of anneal-read cycles was varied between ten and $5{,}120$ in multiples of two, resulting in a variety of idealized and wall-clock QPU execution times, displayed on the horizontal axis.
The annealing time, $\tau$, used in each solution process was also varied from $0.5$ to $312.5$ \SIUnitSymbolMicro s in multiples of five.
The solution qualities and execution times resulting from these annealing time parameterizations are depicted by the five grayscale curves.

It is first apparent that, for a fixed number of anneal-read cycles, as the annealing time is increased, the solution quality also increases.
This is encouraging, as it indicates that the QA hardware is behaving similarly to what the theory of QA predicts.
However, for each increase in annealing time, there are consistent diminishing returns in terms of solution quality.
In particular, between $\tau = 62.5$ \SIUnitSymbolMicro s and $\tau = 312.5$ \SIUnitSymbolMicro s, the differences in solution quality appear mostly inconsequential.
It is also interesting to observe that there are minimal changes in the ideal QPU execution time when using different annealing times for a fixed number of anneal-read cycles.
This is due to an implementation detail of the QPU, where reading solutions takes around $100$ \SIUnitSymbolMicro s, which dominates the ideal execution times presented in Figure \ref{figure:annealing-time-performance} when $\tau$ is less than $100$ \SIUnitSymbolMicro s.
It is also interesting to observe that the wall-clock QPU execution times (i.e., times including service overheads) are much longer than idealized times, typically ranging between two and eight seconds.
These wall-clock results further suggest that increasing the annealing time typically has little effect on the overall realized execution time, as most of the time is dominated by the selected number of anneal-read cycles and service overheads.

Given the diminishing solution quality as annealing time increases, as well as marked execution time increases for $\tau = 312.5$ \SIUnitSymbolMicro s, in this work, we adopt an annealing time of $\tau = 62.5$ \SIUnitSymbolMicro s, which strikes a balance between solution quality and execution time.
This selection is reflected in Section \ref{section:optimization-performance-analysis}.

\section{Additional Run Time Analysis}
\label{section:additional-runtime-analysis}

In this section, we expand on the computational results presented in Section \ref{section:optimization-performance-analysis} by showing the performance results over all of the solution methods considered in this appendix and conducting an average case analysis, accompanying a characteristic example.
These additional results serve to motivate the selection of the five methods presented in Section \ref{section:optimization-performance-analysis}, as well as to show that the characteristic example presented in this work is a suitable representative, not an outlier, of the \ref{eq:cbfm-avg} instance class.

To review the computational setting, all of the classical optimization algorithms were executed on a system with two Intel Xeon E5-2695 v4 processors, each with 18 cores at 2.10 GHz, and 125 GB of memory.
We present an evaluation of the above optimization techniques on problem instances of the largest \ref{eq:cbfm-avg} Ising models that we considered on the \texttt{Advantage\_system4.1} system, containing $5{,}387$ variables.
For each solution technique, parameters that control the execution time of the algorithm (e.g., number of reads for QA, the number of sweeps in SA, or the wall-clock time limit of an IQP method) were varied per the parameterizations described throughout Appendix \ref{section:all-solution-methods} to better understand their effects on solution quality.
The plots show the relative difference in the solution quality (to the best-known solution) over the solve time corresponding to a parameterization.

\begin{figure}[t]
    \includegraphics[width=\textwidth]{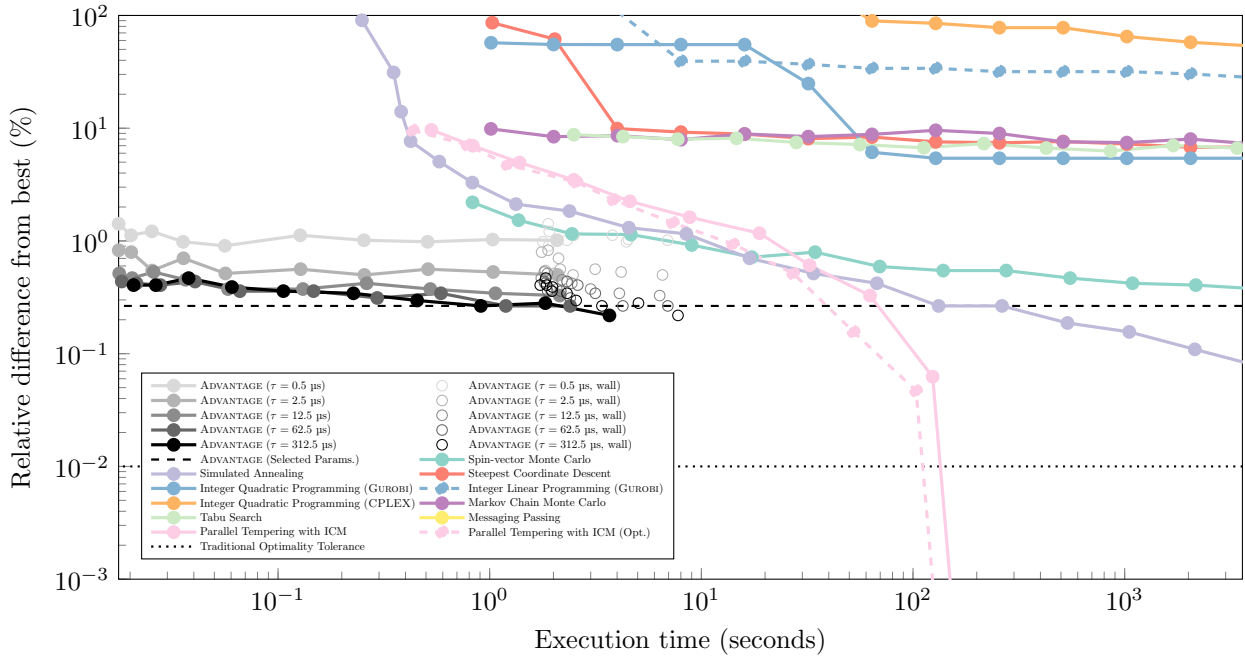}
    \caption{Evaluation of solution quality for a characteristic problem of the \ref{eq:cbfm-avg} instance class with $5{,}387$ decision variables ({\sc Pegasus} lattice size sixteen).}
    \label{figure:characteristic-benchmark-all}
\end{figure}

The benchmarking results for all methods on \ref{eq:cbfm-avg} instance ``16'' (of $50$) are illustrated in Figure \ref{figure:characteristic-benchmark-all}.
A summary of these results follows the discussion of Section \ref{section:optimization-performance-analysis}.
However, this figure justifies our exclusion of some methods.
In particular, we observe that: (1) the \textsc{CPLEX} IQP and \textsc{Gurobi} ILP models are not competitive compared to the \textsc{Gurobi} IQP model and do not bring added insights, and (2) the heuristic methods of tabu search and steepest coordinate descent become stuck in similar local minima, providing limited additional insight.
Message passing provides a minimum difference of $162\%$ and does not appear in the plot's range.
It can also be argued that SVMC brings limited additional insight over the SA heuristic.
However, we elected to include this in Section \ref{section:optimization-performance-analysis}, given the established history of comparing QA to SVMC in the literature.
For heuristic methods, SA and PT-ICM provide the most interesting run time and quality tradeoffs and were often the only methods able to match solutions found by the QA hardware.
Overall, we made the following selections: steepest coordinate descent for representing a greedy approach, \textsc{Gurobi} IQP for representing an off-the-self complete search approach, SVMC by precedent, and SA and PT-ICM as state-of-the-art alternatives.

\begin{figure}[t]
   \includegraphics[width=\textwidth]{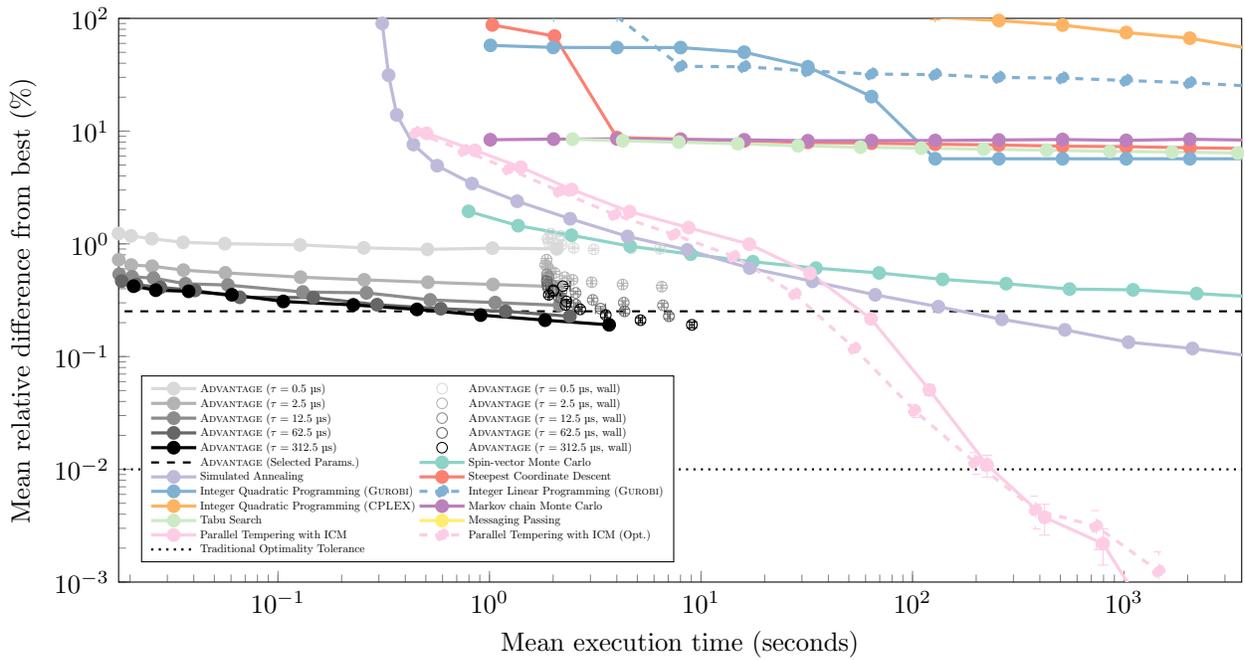}
   \caption{Evaluation of average solution quality over $50$ \ref{eq:cbfm-avg} instances with $5{,}387$ decision variables ({\sc Pegasus} lattice size sixteen).
   Each point corresponds to a mean relative difference, computed over $50$ individual relative differences from the best solution obtained for each instance. Vertical and horizontal error bars correspond to the standard errors of each mean relative difference and each mean solve time, respectively, but tend to be difficult to distinguish because of their small magnitudes.}
   \label{figure:average-benchmark-all}
\end{figure}

Given that \ref{eq:cbfm-avg} defines a distribution of randomly generated optimization instances, it is important to quantify how stable the performance characteristics are across multiple instances from this class.
Figure \ref{figure:average-benchmark-all} presents such an analysis.
It is similar to the other analyses that have been presented but shows performance \emph{averages} across $50$ randomly generated instances.
Error bars are used to show the standard error across multiple \ref{eq:cbfm-avg} instances, but they tend to be difficult to distinguish, as the variance of results across instances is often small.
Message passing again provides a minimum difference of $162\%$ and does not appear in the plot's range.
Overall, these average case results show trends very similar to the characteristic case that is presented in the rest of the paper, providing a strong indication of the stability of the performance characteristics of \ref{eq:cbfm-avg} instances.

\subsection{Solution Quality Details}
\label{subsection:solution-quality-details}

To further elaborate on the consistency of the results presented in Section \ref{section:optimization-performance-analysis} and support the research community in future benchmarking efforts, Table \ref{tbl:obj-values} in this section details the best-known objective values for the $50$ largest \ref{eq:cbfm-avg} instances that were considered in this work.
These instances consist of $5{,}387$ decision variables and are provided as part of the paper's supplementary material for future study.
In particular, Table \ref{tbl:obj-values} reports the {\em best} objective values found by the representative solution methods of Section \ref{section:optimization-performance-analysis} and their parameterizations described in  and Appendix \ref{section:all-solution-methods}.
These detailed results further highlight the consistency in the \ref{eq:cbfm-avg} instances at this problem size, as the qualitative behavior of each solution approach is consistent across all of the instances.
It is important to also note that this table only provides {\em best-known} solutions for these fifty problem instances.
It was observed that the lower bounds produced by the IQP methods were very weak.
Developing tight lower bounds and optimality proofs for these instances remains an important topic for future work.

\begin{table}
    \small
    \centering
    \begin{tabular}{|c|c||c|c|c|c|c|c|}
        \hline
        \bfseries Instance & \bfseries Best Known & \bfseries QA (Avg. 4.1) & \bfseries SA & \bfseries PT-ICM & \bfseries SVMC & \bfseries IQP & \bfseries SCD
        \begin{filecontents*}{data/Pegasus-Lattice_Size-16-table.csv}
instance,best,dw_best,sa_best,pt_best,svmc_best,iqp_best,scd_best
1,-12757.0,-12723.0,-12745.0,-12757.0,-12717.0,-12107.0,-11881.0
2,-12803.0,-12791.0,-12799.0,-12803.0,-12781.0,-12115.0,-11923.0
3,-12763.0,-12717.0,-12745.0,-12763.0,-12709.0,-11933.0,-11917.0
4,-12839.0,-12819.0,-12831.0,-12839.0,-12803.0,-12115.0,-11987.0
5,-12730.0,-12692.0,-12708.0,-12730.0,-12688.0,-12032.0,-11862.0
6,-12828.0,-12796.0,-12810.0,-12828.0,-12782.0,-12024.0,-11936.0
7,-12775.0,-12739.0,-12769.0,-12775.0,-12723.0,-12087.0,-11907.0
8,-12892.0,-12878.0,-12888.0,-12892.0,-12862.0,-12086.0,-11988.0
9,-12699.0,-12665.0,-12679.0,-12697.0,-12645.0,-11977.0,-11821.0
10,-12959.0,-12933.0,-12947.0,-12959.0,-12929.0,-12327.0,-12099.0
11,-12756.0,-12720.0,-12738.0,-12756.0,-12710.0,-11900.0,-11858.0
12,-12736.0,-12714.0,-12722.0,-12736.0,-12712.0,-12110.0,-11870.0
13,-12704.0,-12666.0,-12680.0,-12704.0,-12670.0,-11988.0,-11840.0
14,-12828.0,-12786.0,-12808.0,-12828.0,-12778.0,-12088.0,-11924.0
15,-12882.0,-12856.0,-12878.0,-12882.0,-12838.0,-12208.0,-11972.0
16,-12819.0,-12785.0,-12809.0,-12819.0,-12771.0,-12125.0,-11951.0
17,-12851.0,-12821.0,-12841.0,-12851.0,-12799.0,-12205.0,-12025.0
18,-12834.0,-12802.0,-12818.0,-12834.0,-12804.0,-12214.0,-11956.0
19,-12994.0,-12936.0,-12980.0,-12994.0,-12936.0,-12266.0,-12096.0
20,-12712.0,-12684.0,-12706.0,-12712.0,-12668.0,-12016.0,-11852.0
21,-12854.0,-12838.0,-12844.0,-12854.0,-12818.0,-12096.0,-11956.0
22,-12854.0,-12840.0,-12848.0,-12854.0,-12834.0,-12150.0,-11930.0
23,-12825.0,-12797.0,-12815.0,-12825.0,-12789.0,-12079.0,-11989.0
24,-12820.0,-12780.0,-12802.0,-12820.0,-12770.0,-12056.0,-11956.0
25,-12731.0,-12703.0,-12717.0,-12731.0,-12695.0,-12003.0,-11857.0
26,-12869.0,-12835.0,-12853.0,-12869.0,-12807.0,-12109.0,-11993.0
27,-12908.0,-12884.0,-12900.0,-12908.0,-12872.0,-12126.0,-12062.0
28,-12754.0,-12714.0,-12740.0,-12754.0,-12710.0,-11936.0,-11910.0
29,-12768.0,-12734.0,-12748.0,-12768.0,-12736.0,-11934.0,-11890.0
30,-12778.0,-12754.0,-12772.0,-12778.0,-12726.0,-12026.0,-11866.0
31,-12852.0,-12830.0,-12840.0,-12852.0,-12826.0,-12074.0,-11980.0
32,-12857.0,-12813.0,-12839.0,-12857.0,-12805.0,-12193.0,-11973.0
33,-12848.0,-12820.0,-12834.0,-12848.0,-12798.0,-12040.0,-11966.0
34,-12851.0,-12825.0,-12833.0,-12851.0,-12819.0,-12045.0,-11939.0
35,-12870.0,-12840.0,-12864.0,-12870.0,-12838.0,-12130.0,-12020.0
36,-12798.0,-12788.0,-12794.0,-12798.0,-12770.0,-12158.0,-11932.0
37,-12796.0,-12758.0,-12780.0,-12796.0,-12750.0,-12076.0,-11926.0
38,-12883.0,-12857.0,-12867.0,-12883.0,-12857.0,-12149.0,-11999.0
39,-12849.0,-12837.0,-12841.0,-12849.0,-12827.0,-12217.0,-12003.0
40,-12909.0,-12899.0,-12905.0,-12909.0,-12895.0,-12037.0,-12023.0
41,-12829.0,-12807.0,-12827.0,-12829.0,-12805.0,-12117.0,-11925.0
42,-12710.0,-12692.0,-12700.0,-12710.0,-12690.0,-12036.0,-11834.0
43,-12727.0,-12699.0,-12717.0,-12727.0,-12681.0,-11977.0,-11861.0
44,-12863.0,-12839.0,-12849.0,-12863.0,-12827.0,-12097.0,-12013.0
45,-12822.0,-12788.0,-12814.0,-12822.0,-12772.0,-12166.0,-11936.0
46,-12741.0,-12721.0,-12731.0,-12741.0,-12719.0,-11977.0,-11871.0
47,-12830.0,-12792.0,-12816.0,-12830.0,-12772.0,-12210.0,-11988.0
48,-12877.0,-12847.0,-12873.0,-12877.0,-12825.0,-12137.0,-11997.0
49,-12792.0,-12774.0,-12782.0,-12792.0,-12756.0,-12088.0,-11914.0
50,-12806.0,-12774.0,-12788.0,-12806.0,-12758.0,-12048.0,-11982.0
\end{filecontents*}
\csvreader[head to column names,respect underscore=true,
        late after last line=\ ,
        ]{data/Pegasus-Lattice_Size-16-table.csv}{
            best=\best,dw_best=\dw,sa_best=\sa,pt_best=\pt,svmc_best=\svmc,iqp_best=\iqp,scd_best=\scd}
        {\\\hline\instance & \best & \dw & \sa & \pt & \svmc & \iqp & \scd} \\ \hline
    \end{tabular}
    \caption{A summary of the best objective values found among the representative solution methods on the largest instances considered in this work, with $5{,}387$ decision variables (i.e., {\sc Pegasus} lattice size sixteen).}
    \label{tbl:obj-values}
\end{table}

\end{document}